\documentclass[10pt,draft]{amsart}
\usepackage{amscd}
\usepackage{amssymb}

\newtheorem{thm}[equation]{Theorem}
\newtheorem{cor}[equation]{Corollary}
\newtheorem{prop}[equation]{Proposition}
\newtheorem{lem}[equation]{Lemma}

\theoremstyle{definition}
\newtheorem{dfn}[equation]{Definition}
\newtheorem{rem}[equation]{Remark}
\newtheorem{exa}[equation]{Example}

\numberwithin{equation}{section}

\newcommand{\inj}{\hookrightarrow}

\newcommand{\xar}[1]{\xrightarrow{#1}}
\newcommand{\opn}{\operatorname}
\newcommand{\Cdim}{\opn{Cdim}}

\newcommand{\cat}[1]{\operatorname{\mathsf{#1}}}

\newcommand{\rmitem}[1]{\item[\text{\textup{(#1)}}]}

\newcommand{\mfrak}[1]{\mathfrak{#1}}

\newcommand{\msf}[1]{\mathsf{#1}}

\newcommand{\mrm}[1]{\mathrm{#1}}
\newcommand{\mbb}[1]{\mathbb{#1}}

\newcommand{\tup}[1]{\textup{#1}}

\DeclareMathSymbol{\mathbbk}{\mathord}{AMSb}{"7C}
\renewcommand{\k}{\mathbbk}
\newcommand{\Ext}{\opn{Ext}}
\newcommand{\Hom}{\opn{Hom}}
\newcommand{\fm}{{\mathfrak m}}
\newcommand{\fn}{{\mathfrak n}}
\newcommand{\fp}{{\mathfrak p}}
\newcommand{\fq}{{\mathfrak q}}

\title{Multiplicities of Indecomposable Injectives}
\author{Amnon Yekutieli and James J.\ Zhang }
\address{Department of Mathematics,
Ben Gurion University, Be'er Sheva 84105, ISRAEL}
\email{amyekut@math.bgu.ac.il}
\address{Department of Mathematics, University of
Washington, Seattle, WA 98195, USA}
\email{zhang@math.washington.edu}

\date{6 March, 2004}

\begin{document}

\begin{abstract}
Several results about the multiplicities of indecomposable injectives 
in the minimal injective resolution of a ring exist in the literature.
Mostly these apply to universal enveloping algebras of finite 
dimensional solvable Lie algebras, and to Gorenstein noetherian PI 
local rings. We unify these results and extend them to the much wider
class of rings with Auslander dualizing complexes.
\end{abstract}

\maketitle


\setcounter{section}{-1}
\section{introduction}

Dualizing complexes over noncommutative rings were introduced
by the first author in 1991 \cite{Ye1}. Since then they have been a 
very useful tool in studying noncommutative rings. The structure 
of dualizing complexes is complicated, and is still 
not well understood, except for some special cases. 
In this paper we study a numerical invariant of dualizing 
complexes: the multiplicities of indecomposable 
injectives in a minimal injective resolution of a dualizing 
complex. 

Assume $A$ is a Gorenstein algebra (i.e.\ the bimodule $A$ has 
finite injective dimension on both sides). Then $A$ is a 
dualizing complex over itself, and we obtain information about 
the multiplicities of injectives in a minimal injective resolution 
of $A$. This question has been studied by various authors since 
the 1980's. Some nice results were proved by Barou-Malliavin 
\cite{BM}, Malliavin \cite{Mal}, and Brown-Levasseur \cite{BL} 
for universal enveloping algebras of finite dimensional solvable 
Lie algebras; and by Brown \cite{Br}, Brown-Hajarnavis \cite{BH} 
and Stafford-Zhang \cite{SZ} for noetherian PI algebras. The last 
term of the minimal injective resolution of $A$ was also studied, 
by the above authors and by Ajitabh-Smith-Zhang \cite{ASZ}. The 
theme of this paper is to generalize and unify these existing 
results.

Let $\k$ be a base field. Throughout we shall use the term 
`ring' to refer to a unital associative $\k$-algebra. An 
`$A$-module' shall mean left module, unless stated
otherwise.

Given a left noetherian ring $A$ and a prime ideal 
$\fp \subset A$ let $E(A / \fp)$ denote an injective 
envelope of the $A$-module $A / \fp$, and let $J(\fp)$ denote an 
indecomposable summand of $E(A / \fp)$. So 
$E(A / \fp) \cong J(\fp)^{r}$, 
where $r := \opn{Grank} A / \fp$, the Goldie rank of $A / \fp$. 
Recall that $r$ can be calculated as follows:
the ring of fractions of $A / \fp$ is isomorphic to a matrix ring
$\mrm{M}_r(D)$ over a division ring $D$.

Recall that a {\em minimal injective complex} 
of $A$-modules is a bounded below complex of injectives 
\[ I = (\cdots \to 0 \to I^{i_0} \to \cdots \to  I^i 
\xar{\partial^i} I^{i + 1} \to \cdots) , \] 
such that for every index $i$ the module of cocycles 
$\opn{Ker}(\partial^i)$
is an essential submodule of $I^i$. 
Given a complex $M \in \msf{D}^{+}(\cat{Mod} A)$, 
a {\em minimal injective resolution of $M$} 
is by definition a quasi-isomorphism $M \to I$, where
$I$ is a minimal injective complex. A minimal injective resolution 
of $M$ always exists \cite[4.2]{Ye1}, and it is unique (up to 
a non-unique isomorphism of complexes). For a single module $M$ we 
recover the usual notion of minimal injective resolution.

By `filtered ring' we mean a ring $A$ equipped with an ascending 
filtration \linebreak  $\{ F_i A \}_{i \in \mbb{Z}}$ by 
$\k$-submodules 
such that $F_{-1} A = 0$, $A = \bigcup F_i A$ and
$F_i A \cdot F_j A \subset F_{i + j} A$. The graded ring 
$\opn{gr} A$ is called connected if $\opn{gr}_0 A = \k$
and each $\opn{gr}_i A$ is a finite $\k$-module.

\begin{thm}
\label{xx0.1}
Let $A$ be a filtered ring such that the associated graded ring 
$\opn{gr} A$ is connected, noetherian and has 
an Auslander balanced dualizing 
complex. Let $R$ be the rigid dualizing complex of $A$, and let
$R \to I$ be a minimal injective resolution of $R$ as complex of 
$A$-modules. Then for every prime ideal $\fp \subset A$ 
the multiplicity $\mu_i(\fp)$ of $J(\fp)$ in $I^i$ is
\[ \mu_i(\fp) = 
\begin{cases} 
\opn{Grank} A / \fp & \text{\em if }  
i = - \opn{Cdim}_{R; A} A / \fp , \\
0 & \text{\em otherwise}  .
\end{cases} \]
\end{thm}

The definitions of dualizing complex and related 
terms (such as the Auslander property)
will be reviewed in Section 1. If $R$ is an Auslander dualizing 
complex over $A$ there is a dimension function
$\opn{Cdim}_{R; A} : \cat{Mod} A \to \mbb{Z}$
called the {\em canonical dimension}. For a finite module $M$ the 
formula is 
\[ \opn{Cdim}_{R; A} M := - \opn{inf}\, \{ i \mid 
\opn{Ext}^i_A(M, R) \neq 0 \} . \]

Since the hypothesis of the theorem is left-right symmetric
(replacing $A$ with its opposite ring $A^{\mrm{op}}$)
we also know the multiplicities in the minimal injective 
resolution of $R$ as complex of $A^{\mrm{op}}$-modules.

If $A$ is filtered and the graded ring $\opn{gr} A$ is connected,
noetherian and commutative 
(or somewhat commutative as in Corollary \ref{xx0.2} 
below), then $\opn{gr} A$ has an Auslander balanced dualizing 
complex, and $\opn{Cdim} = \opn{GKdim}$ \cite[6.9]{YZ1}. Therefore 
we have the following result.

\begin{cor} \label{xx0.2} 
Let $A$ be a filtered ring such that the associated 
graded ring $\opn{gr} A$ is connected, noetherian and PI 
\tup{(}or FBN, or has enough normal elements in the sense of 
\cite[p.\ 36]{YZ1}\tup{)}. 
Let $R$ be a rigid dualizing complex over $A$, and let
$R \to I$ be a minimal injective resolution of $R$ as complex of 
$A$-modules.
Then for every prime ideal $\fp \subset A$ 
the multiplicity $\mu_i(\fp)$ of $J(\fp)$ in $I^i$ is
\[ \mu_i(\fp) =
\begin{cases} 
\opn{Grank} A / \fp & \text{\em if }   
i = - \opn{GKdim} A / \fp , \\
0 & \text{\em otherwise} .
\end{cases} \]
\end{cor}

Suppose $A$ is the universal enveloping algebra $U(\mfrak{g})$ 
of a finite dimensional Lie algebra $\mfrak{g}$. 
Then the associated graded 
ring $\opn{gr} A$ is commutative, and the rigid dualizing 
complex of $A$ is $A^\sigma[d]$ for some automorphism 
$\sigma$, where $d : = \opn{dim}_{\k} \mfrak{g}$. 
Hence:

\begin{cor} \label{xx0.3} 
Let $A := U(\mfrak{g})$, the universal enveloping algebra of 
a Lie algebra $\mfrak{g}$ with 
$\opn{dim}_{\k} \mfrak{g} = d < \infty$. 
Let $A \to I$ be a minimal injective resolution of $A$ as 
$A$-module. Then for every prime ideal $\fp \subset A$ 
the multiplicity $\nu_i(\fp)$ of $J(\fp)$ in $I^i$ is
\[ \nu_i(\fp) =
\begin{cases} 
\opn{Grank} A / \fp & \text{\em if} \quad 
i = d - \opn{GKdim} A / \fp , \\
0 & \text{\em otherwise} .
\end{cases} \]
\end{cor}

Corollary \ref{xx0.3} extends a result of Malliavin 
\cite[4.17]{Mal} and Brown-Levasseur \cite[6.1]{BL}, which was 
proved for enveloping algebras of solvable Lie algebras. 
Note that our result is valid for any Lie algebra. Corollary 
\ref{xx0.2} also applies to factor rings of $U(\mfrak{g})$ and many 
quantized algebras listed in \cite{GL}. 

The universal enveloping algebra of a finite dimensional Lie 
algebra is Auslander regular of finite global dimension. Most 
algebras in this paper do not have finite global (or injective) 
dimension, but often admit Auslander dualizing complexes. Such 
algebras are generalization of Auslander Gorenstein or Auslander 
regular algebras. Basic properties of Auslander dualizing 
complexes and their applications are given in \cite{Ye2}, 
\cite{YZ1}, \cite{WZ1}.

For FBN rings, every indecomposable injective is of the form 
$J(\fp)$ for some prime ideal $\fp$. In this case we can say more 
about the structure of dualizing complexes. 
A version of Theorem \ref{xx0.4} for non-FBN 
rings is given in Proposition \ref{xx3.1}.

\begin{thm}
\label{xx0.4}
Let $A$ and $B$ be two FBN rings, and let $R$ be an Auslander, 
weakly bifinite, $\opn{Cdim}$-symmetric dualizing complex 
over $(A,B)$. Let $R \to I$ be a minimal injective 
resolution of $R$ as complex of $A$-modules. Then for every $i$ we 
have
\[ I^i \cong \bigoplus\, J(\fp)^{\mu_i(\fp)} , \]
where the sum ranges over all prime ideals $\fp \subset A$ with 
$\opn{Cdim}_{R; A} A / \fp = -i$. 
The multiplicity $\mu_i(\fp)$ is a \tup{(}finite\tup{)}
positive integer.
\end{thm}

Theorem 0.4 is a generalization of the results of Brown 
\cite[Theorem C]{Br}, and Stafford-Zhang \cite[3.15]{SZ}. 
This theorem is part of \cite[4.10]{YZ2}. However, the 
statement presented here is implicit in \cite[4.10]{YZ2}
and the proof we will give is slightly different from
the one in \cite[4.10]{YZ2}. The 
structure of dualizing complexes over FBN rings was studied 
further in \cite{YZ2} (see Remark \ref{xx5.4}). In some special 
cases, we are able to compute the multiplicity $\mu_i(\fp)$ in 
Theorem \ref{xx0.4}. The following corollary is analogous to 
\cite[5.5]{BH}.

\begin{cor}
\label{xx0.5} Let $A$ be a noetherian, complete, local, PI 
$\k$-algebra with injective dimension $d < \infty$ and maximal 
ideal $\fm$. Suppose that $A / \fm$ is finite over $\k$.
Let $A \to I$ be a minimal injective 
resolution of $A$ as $A$-module.
Then for every $i$ we have an isomorphism
\[ I^i \cong \bigoplus\, J(\fp)^{\opn{Grank} A / \fp} , \]
where the sum ranges over all prime ideals $\fp \subset A$ with 
$\opn{Kdim} A / \fp = d - i$. 
\end{cor}

Theorem \ref{xx0.4} is proved in Section 3, and the 
other theorems and corollaries are proved in Section 5.

\medskip \noindent 
{\bf Acknowledgments.} Both authors were supported by the 
US-Israel Binational Science Foundation. The second author was 
also supported by the NSF, a Sloan Research Fellowship and the 
Royalty Research Fund of the University of Washington. Authors
also thank the referee for careful reading of the manuscript 
and valuable comments.


\section{Review of Definitions}

Let $A$ be a ring (i.e.\ a $\k$-algebra). As mentioned earlier, 
we shall usually work with left $A$-modules. We treat right 
$A$-modules as $A^{\mrm{op}}$-modules, where $A^{\mrm{op}}$ is 
the opposite ring. Bimodules are modules over the enveloping ring 
$A^{\mrm{e}} := A \otimes A^{\mrm{op}}$, where 
$\otimes := \otimes_{\k}$. 
A finitely generated $A$-module is simply called 
finite. We refer to \cite{VdB,Ye1,Ye2,YZ1} for basic material on 
complexes and derived categories, and to \cite{GW, MR} for basic 
ring theoretic notions such as PI and FBN. 

The category of $A$-modules is denoted by 
$\cat{Mod} A$, and its derived category is 
$\msf{D}(\cat{Mod} A)$.

We now recall several notions related to dualizing complexes.

\begin{dfn} \label{xx1.1} \cite[3.3]{Ye1}, \cite[1.1]{YZ1}
Assume $A$ is a left noetherian $\k$-algebra and $B$ is a right 
noetherian $\k$-algebra. A complex 
$R \in \cat{D}^{\mrm{b}}(\cat{Mod} A \otimes B^{\mrm{op}})$ 
is called a {\em dualizing complex} over $(A,B)$ 
if it satisfies the following conditions:
\begin{enumerate}
\rmitem{i} $R$ has finite injective dimension over $A$ and 
$B^{\mrm{op}}$.
\rmitem{ii} $R$ has finite cohomology modules
over $A$ and $B^{\mrm{op}}$.
\rmitem{iii} The canonical morphisms $B \to \opn{RHom}_A(R, R)$ in
$\cat{D}(\cat{Mod} B^{\mrm{e}})$, and 
$A \to \opn{RHom}_{B^{\mrm{op}}}(R, R)$ in 
$\cat{D}(\cat{Mod} A^{\mrm{e}})$, are both isomorphisms.
\end{enumerate}
 
If moreover $A = B$ we say $R$ is a dualizing complex over $A$.
\end{dfn}

Whenever we say $R$ is a dualizing complex over $(A,B)$ we 
tacitly assume that $A$ is left noetherian and $B$ is right 
noetherian.

Let $R$ be a dualizing complex over $(A,B)$ and let $M$ be a 
finite $A$-module. The {\em grade} of $M$ with respect to 
$R$ is defined to be
\[ j_R(M) = \inf\, \{q \mid \opn{Ext}^q_A(M, R) \neq 0 \} . \]
The grade of a $B^{\mrm{op}}$-module is defined similarly.

\begin{dfn} \label{xx1.2}\cite[1.2]{Ye2}, \cite[2.1]{YZ1}
A dualizing complex $R$ over $(A,B)$ is called
{\em Auslander} if
\begin{enumerate}
\rmitem{i} For every finite $A$-module $M$, every integer
$q$ and every $B^{\mrm{op}}$-submodule 
$N \subset \opn{Ext}^q_A(M, R)$ one has 
$j_R(N) \geq q$.
\rmitem{ii} The same holds after exchanging $A$ and 
$B^{\mrm{op}}$. 
\end{enumerate}
\end{dfn}

The canonical dimension of a finite $A$-module $M$ with respect 
to the Auslander dualizing complex $R$ is defined to be
\[ \opn{Cdim}_{R; A} M := - j_R(M) . \]
By \cite[2.10]{YZ1}, $\opn{Cdim}_{R; A}$ 
is a finitely partitive, exact dimension function
(cf.\ \cite[6.8.4]{MR}). 
By the left-right symmetry of the situation we obtain a dimension 
function $\opn{Cdim}_{R; B^{\mrm{op}}}$ on
$\cat{Mod} B^{\mrm{op}}$. When there is no danger of confusion 
we will tend to drop some of the subscripts from the expressions  
$\opn{Cdim}_{R; A}$  and $\opn{Cdim}_{R; B^{\mrm{op}}}$.

\begin{dfn} \label{xx1.3}
A dualizing complex $R$ over $(A,B)$ is called {\it bifinite} 
(respectively {\it weakly bifinite}) if the following conditions 
hold:
\begin{enumerate}
\rmitem{i} For every $A$-bimodule $M$ 
which is finite on both sides 
(respectively, and is a subquotient of $A$) the 
$(A, B)$-bimodule $\opn{Ext}^q_A(M, R)$ is finite on both sides.
\rmitem{ii} The same holds after exchanging $A$ and $B^{\mrm{op}}$.
\end{enumerate}
\end{dfn}

Not every dualizing complex is weakly bifinite (see Example 
\ref{xx6.1}).

The notion of rigidity is due to Van den Bergh \cite{VdB}.

\begin{dfn} \label{xx1.4}\cite[8.1]{VdB}
Suppose $R$ is a dualizing complex over $A$. If there is an 
isomorphism
\[ \rho : R \to \opn{RHom}_{A^{\mrm{e}}}(A, R \otimes R) \]
in $\cat{D}(\cat{Mod} A^{\mrm{e}})$ then we call $(R,\rho)$, 
or simply $R$, a {\it rigid} dualizing complex.
\end{dfn}

\begin{dfn} \label{xx1.5} \cite[4.5]{Ye1}
A dualizing complex $R$ over $(A,B)$ is called 
{\it pre-balanced} if
\begin{enumerate}
\rmitem{i} For every simple $A$-module $M$ one has
$\opn{Ext}^i_A(M, R) = 0$ for all $i \neq 0$,
and $\Ext^0_A(M, R)$ is a simple $B^{\mrm{op}}$-module. 
\rmitem{ii} The same holds after 
exchanging $A$ and $B^{\mrm{op}}$.
\end{enumerate}
\end{dfn}

\begin{rem} \label{xx1.6} 
Pre-balanced dualizing complexes often exist for semilocal or FBN 
rings. But this is not true in general. For example, the $n$th 
Weyl algebra, for $n \geq 2$, does not admit any pre-balanced 
dualizing complexes (see Example \ref{xx6.5}), though it 
naturally admits an Auslander rigid dualizing complex. 
\end{rem}

When $A$ is connected graded, there is a notion of balanced
dualizing complex introduced in \cite{Ye1}. Let $\fm=A_{>0}$ and
let $A'$ be the graded vector space dual of $A$.
Let $\Gamma_{\fm}$ denote the torsion functor $\lim_{n\to \infty} 
\opn{Hom}_A(A/\fm^n, -)$. 

\begin{dfn}\cite[4.1]{Ye1} \label{xx1.7} 
A dualizing complex $R$ over a connected graded ring $A$ is called 
{\em balanced} if there are isomorphisms 
$$\opn{R\Gamma}_{\fm}(R)\cong A'
\cong \opn{R\Gamma}_{\fm^\circ}(R)$$
in $\cat{D}^{\mrm{b}}(\cat{Mod} A^e)$. A balanced dualizing complex
over a local ring is defined similarly \cite[Definition 3.7]{CWZ}.
\end{dfn}

\begin{dfn} \label{xx1.8} 
Suppose we are given a dimension function $\opn{dim}_A$ 
on $\cat{Mod} A$ and a dimension function 
$\opn{dim}_{B^{\mrm{op}}}$ on $\cat{Mod} B^{\mrm{op}}$. The pair
$(\opn{dim}_A, \opn{dim}_{B^{\mrm{op}}})$ is denoted by 
$\opn{dim}$. 
\begin{enumerate}
\rmitem{i} We say that $\opn{dim}$ is {\em symmetric} if
for every $(A,B)$-bimodule $M$ finite on both sides there
is equality
\[ \opn{dim}_A M = \opn{dim}_{B^{\mrm{op}}} M . \]
\rmitem{ii} In case $A = B$ we say that 
$\opn{dim}$ is {\em weakly symmetric} if 
for every $A$-bimodule $M$ that is a subquotient of $A$ there
is equality
\[ \opn{dim}_A M = \opn{dim}_{A^{\mrm{op}}} M . \]
\rmitem{iii} Assume $R$ is an Auslander dualizing complex
over $(A,B)$ (respectively over $A$). 
We say  $R$ is {\it $\opn{Cdim}$-symmetric} (respectively 
{\it $\opn{Cdim}$-weakly symmetric}) if the pair 
$\opn{Cdim}_R := 
(\opn{Cdim}_{R; A}, \opn{Cdim}_{R; B^{\mrm{op}}})$ 
is symmetric (respectively weakly symmetric).
\end{enumerate}
\end{dfn}

The Gelfand-Kirillov dimension $\opn{GKdim}$ is always symmetric  
\cite[8.3.14(ii)]{MR}, and Krull dimension $\opn{Kdim}$ is 
symmetric for FBN rings \cite[6.4.13]{MR}. It follows from 
\cite[4.8]{VdB} and \cite[6.22]{YZ1} that rigid Auslander 
dualizing complexes over filtered algebras are 
$\opn{Cdim}$-weakly symmetric. A similar statement holds in the 
semilocal case \cite[0.1]{WZ1}. But Example \ref{xx6.2} shows 
that not every Auslander bifinite dualizing complex is 
$\opn{Cdim}$-weakly symmetric.


\section{Preliminary results}

In this section we collect some preliminary results about
the pre-balanced condition, bifiniteness and so on. 

\begin{lem} \label{xx2.1} 
The following dualizing complexes are bifinite or weakly bifinite.
\begin{enumerate}
\item If $A$ and $B$ are noetherian PI algebras, 
then every dualizing complex $R$ over $(A,B)$ is bifinite.
\item If $A$ is filtered and $\opn{gr} A$ is noetherian with a
balanced dualizing complex, then the rigid dualizing complex over 
$A$ is weakly bifinite.
\item{} Let $(A,\fm)$ and $(B,\fn)$ be noetherian complete 
semilocal algebras such that there is a Morita duality between
them. Suppose that
\begin{enumerate}
\item $\Gamma_{\fm}$ and $\Gamma_{\fn^\circ}$ have finite cohomogical
dimension,
\item $A$ and $B$ satisfy the $\chi$-condition, and 
\item $A/\fm$ is PI.
\end{enumerate}
Then every pre-balanced dualizing complex over $(A,B)$ is bifinite.
\end{enumerate}
\end{lem}

\begin{proof} 
(1) Let $M$ be an $A$-bimodule finite on both sides. It follows 
from the properties of dualizing complexes that 
$\mrm{Ext}^i_A(M,R)$ is a finite $B^{\mrm{op}}$-module. It remains 
to show that it is finite over $A$. So we may forget the 
$B^{\mrm{op}}$-module structure of $R$ and by the long exact 
sequence we may replace $R$ by a finite $A$-module. Then the 
result follows from \cite[3.5]{SZ}. 

\medskip \noindent
(2) This follows from \cite[6.21]{YZ1}. (However we don't know if 
$R$ is bifinite in this case.) 

\medskip \noindent
(3) This follows from \cite[0.1(4)]{CWZ}. For the definitions 
and some details of the conditions listed, see \cite{CWZ,WZ1}. 
\end{proof}

The pre-balanced condition is a natural one according to 
\cite{CWZ}. From the next lemma we see that it is automatic for 
local rings.  

\begin{lem} \label{xx2.2}
Let $R$ be a weakly bifinite dualizing complex over $(A,B)$. 
Suppose that for every simple $B^{\mrm{op}}$-module $S$ the ring
$B / \opn{Ann}_{B^{\mrm{op}}} (S)$ is simple artinian.
\begin{enumerate}
\item Condition \tup{\ref{xx1.5}(i)} 
implies \tup{\ref{xx1.5}(ii)}. Consequently, $R$ is pre-balanced
if and only \tup{\ref{xx1.5}(i)} holds.
\item If $A$ is local then $B$ is local too, and 
$R[d]$ is pre-balanced for some integer $d$.
\end{enumerate}
\end{lem}

\begin{proof} 
(1) Suppose that \ref{xx1.5}(i) holds. 
 Let $S$ be a simple $B^{\mrm{op}}$-module and 
$\fn := \opn{Ann}_{B^{\mrm{op}}} (S)$. 
By hypothesis, $B / \fn$ is simple artinian, and 
hence $B / \fn \cong S^{r}$ as $B^{\mrm{op}}$-modules. Since 
$R$ is weakly bifinite, 
$ M := \Ext^j_{B^{\mrm{op}}}(B / \fn, R)$ is 
noetherian on both sides. It is a right $(B / \fn)$-module, 
so it's right artinian. By Lenagan's lemma \cite[7.10]{GW}, 
$M$ is also left artinian. It follows that
$A / \opn{Ann}_A (M)$ is an $A$-submodule of a 
finite direct sum of copies of $M$, which is left artinian. By 
\ref{xx1.5}(i),
\[ \Ext^i_A(\Ext^j_{B^{\mrm{op}}}(B / \fn, R),R)
= \opn{Ext}^i_A(M, R) = 0 \]
for all $i \neq 0$. Then the double-Ext spectral sequence
\cite[1.7]{YZ1} shows that 
$\Ext^j_{B^{\mrm{op}}}(B / \fn, R) = 0$ 
for all $j \neq 0$ and that 
\[ B / \fn \cong \Ext^0_A(\Ext^0_{B^{\mrm{op}}}(B / \fn, R), R)  
. \]
Hence 
$\Ext^j_{B^{\mrm{op}}}(S, R) = 0$ 
for all $j \neq 0$ and 
\[ S \cong \Ext^0_A(\Ext^0_{B^{\mrm{op}}}(S, R), R) . \]
It follows from \ref{xx1.5}(i) that the length of the $A$-module
$\Ext^0_{B^{\mrm{op}}}(S, R)$ is
equal to the length of the $B^{\mrm{op}}$-module
$\Ext^0_A(\Ext^0_{B^{\mrm{op}}}(S, R), R)$,
which is one. Therefore $\Ext^0_{B^{\mrm{op}}}(S, R)$ is simple.

\medskip \noindent
(2)  Let $\fm$ be the maximal ideal of $A$. Shift $R$ so that 
$\Ext^0_{A}(A / \fm, R) \neq 0$. By weak finiteness and Lenagan's 
lemma as in the proof of (1), 
$\Ext^i_A(A / \fm, R)$ is artinian on 
both sides for all $i$. Let
\[ i_0 := \min\, \{ i \mid \Ext^i_{A}(A / \fm, R)\neq 0 \}, \quad  
i_1 := \max\, \{ i \mid \Ext^i_{A}(A / \fm, R)\neq 0 \} . \]
By induction on length one sees that 
for all finite length $A$-modules $M$
\[ i_0 = \min\, \{ i \mid \Ext^i_{A}(M, R) \neq 0 \}, \quad 
i_1 = \max\, \{ i \mid \Ext^i_{A}(M, R) \neq 0 \} . \]

Let $S$ be a simple $B^{\mrm{op}}$-module. By our hypothesis
$\fn := \opn{Ann}_{B^{\mrm{op}}} (S)$ is an ideal 
such that $B / \fn$ is simple 
artinian. For each $j$ let 
$M_j := \Ext^j_{B^{\mrm{op}}}(B / \fn, R)$. 
Then by weak finiteness and Lenagan's lemma $M_j$ is 
artinian on both sides. Let 
\[ j_0 := \min\, \{ j \mid M_j \neq 0 \}, \quad 
j_1 := \max\, \{ j \mid M_j \neq 0 \} . \]
Then $\Ext^i_A(M_j, R) \neq 0$ for $i = i_0, i_1$ and 
$j = j_0, j_1$.
Note that these are four corner terms in the $E_2$-page of
the double-Ext spectral sequence \cite[1.7]{YZ1}
\[ E^{ij}_2 := \Ext^i_A(\Ext^j_{B^{\mrm{op}}}(B / \fn, R), R)
\Rightarrow B / \fn . \]
Since the terms $E^{i_0j_1}_2$ and $E^{i_1j_0}_2$ will survive
in the $E_{\infty}$-page, these two terms must be on the diagonal,
namely $i_0 = j_1$ and $i_1 = j_0$. This shows that 
$i_0 = i_1 = j_0 = j_1$,
and all are equal to $0$ because $i_0 \leq 0 \leq i_1$.
As a consequence $\Ext^j_{B^{\mrm{op}}}(S, R) = 0$ 
for all $j \neq 0$,
and $\Ext^0_{B^{\mrm{op}}}(S, R)$ is a nonzero artinian 
$A$-module. Therefore $\Ext^0_A(-, R)$ and 
$\Ext^0_{B^{\mrm{op}}}(-, R)$ are a duality
between the category of finite length
$A$-modules and that of finite length $B^{\mrm{op}}$-modules. 
Therefore $B$ is local and $R$ is pre-balanced.
\end{proof}

Recall from \cite{WZ2} that $A$ satisfies the {\it Nakayama 
condition} if for every nonzero finite right $A$-module $N$
there is a (simple) left $A$-module $S$ such that 
$N \otimes_A S \neq 0$. It is easily seen that noetherian 
semilocal rings and FBN rings satisfy the Nakayama condition. 
One can also verify that the first Weyl algebra, 
which is neither FBN nor semilocal, satisfies 
the Nakayama condition.

The following Proposition is an ungraded version of 
\cite[3.9]{Ye1}.

\begin{prop} \label{xx2.3} 
Let $R_1$ and $R_2$ be pre-balanced dualizing complexes 
over $(A,B)$. If $A$ satisfies the Nakayama condition, then 
there is an invertible $A$-bimodule $L$ such that 
$R_1 \cong L \otimes_{A} R_2$.
\end{prop}

To prove this we need to use two-sided tilting complexes (see 
\cite{Ri1,Ri2,Ye3}). A two-sided tilting complex 
$T \in \cat{D}^{\mrm{b}}(\cat{Mod} B \otimes A^{\mrm{op}})$ 
is called {\it pre-balanced} if for every simple 
$A$-module $S$ and every $i \neq 0$ one has 
$\mrm{H}^{i}(T \otimes^{\mrm{L}}_{A} S)= 0$.

\begin{lem} \label{xx2.4}
If $A$ satisfies the Nakayama condition, then a two-sided 
tilting complex 
$T \in \cat{D}^{\mrm{b}}(\cat{Mod} B \otimes A^{\mrm{op}})$ 
is pre-balanced if and 
only if it is isomorphic to an invertible bimodule. 
\end{lem}

\begin{proof} 
If $T$ is isomorphic to an invertible bimodule, then it is 
clearly pre-balanced. So we only need to prove the other 
implication. 

Let 
$i_{1} := \max\, \{i \mid \mrm{H}^{i} T \neq 0 \}$.
Since $\mrm{H}^{i_{1}} T$ is a finite $A^{\mrm{op}}$-module, 
by the Nakayama condition there is a simple $A$-module 
$S$ such that  
$(\mrm{H}^{i_{1}} T) \otimes_{A} S \neq 0$. 
By the K\"{u}nneth spectral sequence (cf.\ \cite[2.1]{Ye3}) there 
is an isomorphism
\[ (\mrm{H}^{i_{1}} T) \otimes_{A} S \cong
\mrm{H}^{i_{1}}(T \otimes^{\mrm{L}}_{A} S) . \]
Since $T$ is pre-balanced we see that $i_{1} = 0$.

We conclude that $(\mrm{H}^{0} T) \otimes_{A} S \neq 0$
for all simple $A$-modules $S$. Therefore also
$(\mrm{H}^{0} T) \otimes_{A} M \neq 0$
for every nonzero finite $A$-module $M$.

Proceeding exactly as in \cite[2.3]{Ye3} we conclude that 
$T \cong \mrm{H}^{0} T$ is an invertible bimodule.
\end{proof}

\begin{proof}[Proof of Proposition \tup{\ref{xx2.3}}]
First we note that since for $i = 1, 2$ the complex $R_{i}$ 
is pre-balanced, the functor
$M \mapsto \opn{Ext}^{0}_A(M, R_{i})$
is a duality between the category  of finite length $A$-modules
and that of finite length $B^{\mrm{op}}$-modules. 
Furthermore, for every finite length $A$-module $M$, 
$\opn{RHom}_{A}(M, R_{i}) \cong 
\opn{Ext}^{0}_{A}(M, R_{i})$;
and likewise for $B^{\mrm{op}}$-modules.

By \cite[3.9 and 3.10]{Ye1} there is a two-sided
tilting complex $T$ over $A$ such that 
$R_{2} \cong T \otimes^{\mrm{L}}_{A} R_{1}$, 
and such that 
\[ \opn{RHom}_{B^{\mrm{op}}}(\opn{RHom}_{A}(M, R_{1}), R_{2}))
\cong T \otimes^{\mrm{L}}_{A} M \]
for all $M \in \msf{D}^{\mrm{b}}_{\mrm{f}}(\cat{Mod} A)$. 
Taking $M$ to be a simple $A$-module we get
$\mrm{H}^{i}(T \otimes^{\mrm{L}}_{A} M) = 0$
if $i \neq 0$. Thus $T$ is pre-balanced, and by Lemma \ref{xx2.4} 
it is isomorphic to an invertible bimodule. 
\end{proof} 

Let $\fp$ be a prime ideal of $A$. We denote by $J(\fp)$  the 
indecomposable injective module associated to $\fp$, 
i.e.\ the $A$-injective envelope of a uniform left ideal 
of $A / \fp$. 

Given an $(A / \fp)$-module $M$ its reduced rank (at 
the prime $\fp$) is 
\[ \opn{red{.}rank} M 
:= \opn{Grank} (M / \tau(M)) , \] 
where $\tau(M)$ denotes the torsion submodule of $M$ and
$\opn{Grank}$ denotes Goldie rank. If $Q$ is the ring of fractions 
of $A / \fp$ then 
$\opn{red{.}rank} M$ equals the length of 
$Q \otimes_A M$ as $Q$-module. 

A complex $I$ of $A$-modules is called a {\em minimal injective 
complex} if (a) each term $I^i$ is $A$-injective and (b) each $\ker 
\partial^i$ is an essential submodule of $I^i$ where $\partial^i:
I^i\to I^{i+1}$ is the $i$th coboundary map. 
The next lemma is \cite[4.12]{YZ2}, which is a generalization of 
\cite[2.3]{Br}.

\begin{lem} \cite[4.12]{YZ2} \label{xx2.5}
Suppose $A$ is a left noetherian ring. Let $I$ be a minimal 
injective complex of $A$-modules. Let $\mfrak{p}$ be a prime 
ideal of $A$, and let $\mu_{i}(\mfrak{p})$ be the multiplicity of 
$J(\mfrak{p})$ in $I^i$. Then
\begin{enumerate}
\item The image of the map 
\[ \opn{Hom}_A(A/\fp,\partial^{i - 1}): \;\;
\opn{Hom}_{A}(A / \mfrak{p}, I^{i-1}) \to 
\opn{Hom}_{A}(A / \mfrak{p}, I^{i}) \]
is a torsion $(A/\mfrak{p})$-module.
\item There is equality
\[ \mu_{i}(\mfrak{p}) = 
\opn{red{.}rank} \opn{Hom}_{A}(A / \mfrak{p}, I^{i}) =
\opn{red{.}rank} \opn{Ext}^{i}_{A}(A / \mfrak{p}, I) . 
\]
\item Let $\mfrak{p}$ and $\mfrak{q}$ be two primes of $A$ 
and $M$ an $(A / \mfrak{p})$-$(A / \mfrak{q})$-bimodule. Assume $M$ 
is nonzero, torsion-free as $(A / \mfrak{q})^{\mrm{op}}$-module
and finite non-torsion as $(A / \mfrak{p})$-module. If 
$\mu_{i}(\mfrak{p}) \geq 1$ then $\opn{Ext}^{i}_{A}(M, I)$ 
is a non-torsion $(A / \mfrak{q})$-module.
\end{enumerate}
\end{lem}


\section{Multiplicity, Part I}

In this section $R$ is a dualizing complex over a pair of rings
$(A,B)$, and $R \to I$ is a minimal injective resolution of 
$R$ as complex of $A$-modules. For a prime ideal 
$\mfrak{p} \subset A$ and an integer $i$ let $\mu_i(\fp)$ be 
the multiplicity of $J(\fp)$ in $I^i$. 
Using Lemma \ref{xx2.5} we are able to compute the multiplicities
$\mu_i(\fp)$ in some cases. 

\begin{prop} \label{xx3.1}
Suppose $R$ is a weakly bifinite dualizing complex over $(A,B)$. 
Then $\mu_i(\fp)$ is finite for all $i$ and all primes 
$\fp\subset A$. If moreover $R$ is Auslander and 
$\opn{Cdim}$-symmetric, then $\mu_i(\fp)$ is nonzero if and only 
if $i = - \Cdim_{R; A} A / \fp$. 
\end{prop}

\begin{proof} 
Because $R$ is weakly bifinite, $\Ext^i_A(A / \fp, R)$ is a finite  
$A / \fp$ -module. Therefore 
$\opn{red{.}rank}_{A / \mfrak{p}} \Ext^i_A(A / \fp, R)$ 
is finite. The first assertion follows from Lemma \ref{xx2.5}(2).

To prove the second assertion let 
$j := j_R(A / \fp) = -\opn{Cdim}_{R; A} A / \fp$.
By \cite[1.11]{YZ1} we know that
$\mu_i(\fp)\neq 0$ for some $i$. So it suffices 
to show that $\mu_i(\fp)=0$ for all $i\neq j$. But 
when $i \neq j$ we have
\[ \opn{Cdim}_{R; B^{\mrm{op}}} \Ext^i_A(A / \fp, R) <
\opn{Cdim}_{R; B^{\mrm{op}}} \Ext^j_A(A / \fp, R) = -j . \]
Hence $\Ext^i_A(A / \fp, R)$ is a torsion 
$(A / \fp)$-module. The assertion now follows from 
Lemma \ref{xx2.5}(2).
\end{proof}

\begin{proof}[Proof of Theorem \tup{\ref{xx0.4}}]
First of all since $A$ is noetherian, every injective module is a 
direct sum of indecomposable injectives. Since $A$ is FBN, every 
indecomposable injective is of form $J(\fp)$ for some prime ideal 
$\fp \subset A$. 
By Proposition \ref{xx3.1}, $\mu_i(\fp)$ is finite, and it is nonzero 
if and only if $i = -\opn{Cdim}_{R; A} A / \fp$.
\end{proof}

\cite[3.15]{SZ} is an almost immediate consequence:

\begin{cor} \cite[3.15]{SZ} \label{xx3.2} 
Suppose $A$ is a noetherian, AS-Gorenstein, PI algebra. Then the 
$i$th term of the minimal injective resolution of $A$ 
as complex of $A$-modules is isomorphic to
\[ \bigoplus J(\fp)^{\nu_i(\fp)} , \]
where the sum ranges over all prime ideals $\fp$ such that
$\opn{Kdim} A / \fp = \opn{Kdim} A - i$. 
The multiplicity $\nu_i(\fp)$
is a \tup{(finite)} positive integer.
\end{cor}

\begin{proof} 
Since $A$ is AS-Gorenstein and PI it follows from
\cite[3.10]{SZ} that $A$ is 
Auslander Gorenstein and $\opn{Kdim}$-Macaulay. Let $d$ be the 
injective dimension of $A$, which equals its Krull dimension. 
Then $R := A[d]$ is an Auslander, bifinite, 
$\opn{Cdim}$-symmetric dualizing complex over $A$. Now use
Proposition \ref{xx3.1}, noting that 
$\nu_i(\fp) = \mu_{i - d}(\fp)$.
\end{proof}

If $A$ is prime the multiplicities $\mu_i(\mfrak{o})$ for the zero 
ideal $\mfrak{o}$ can be computed, as we now show. 
But first a bit of notation. Say 
$M$ is an $A$-$B$-bimodule. We write 
$\opn{Grank}_A M$ for the Goldie rank of $M$ considered as 
$A$-module, and likewise $\opn{Grank}_{B^{\mrm{op}}} M$.

\begin{prop} \label{xx3.3}
Suppose $A$ and $B$ are prime rings and 
$R$ is an Auslander dualizing complex over $(A,B)$. 
Let $H^{-d} R$ be the lowest nonzero cohomology of $R$. 
\begin{enumerate}
\item Let $Q(A)$ be the ring of fractions of $A$, and likewise
$Q(B)$. Then 
\[ Q(A) \otimes_{A} H^{i} R \cong H^{i} R \otimes_{B} Q(B) \]
as $A$-$B$-bimodules for all $i$, and they are $0$ for 
$i \neq -d$.
\item The bimodule $H^{-d}(R)$ is torsion-free on both sides,
\[ \opn{Grank}_A H^{-d} R =
\opn{Grank}_{B^{\mrm{op}}} B \]
and
\[ \opn{Grank}_{B^{\mrm{op}}} H^{-d} R =
\opn{Grank}_A A . \]
\item  Let $\mfrak{o}$ be the zero ideal of $A$. Then 
\[ \mu_{-d}(\mfrak{o}) = \opn{Grank}_{B^{\mrm{op}}} B \]
and 
$\mu_i(\mfrak{o}) = 0$ for all $i \neq -d$.
\end{enumerate}
\end{prop}

\begin{proof} 
(1) We may assume that $R^{-d-1} = 0$. Hence
$\Ext^{-d}_A(-, R) = \Hom_A(-, H)$ where 
$H := \mrm{H}^{-d}(R)$. Since $A$ and $B$ are prime of 
$\opn{Cdim} = d$, every torsion module $M$ has 
$\opn{Cdim} M < d$. 
Thus, if $M$ is a torsion finite $A$-module, then $\Hom_A(M, H) = 0$. 
Hence $H$ is torsion-free on the left; and by symmetry also on the 
right. Let $b \in B$ be a regular element. 
Then right multiplication by $b$ 
on $Q(A) \otimes_{A} H$ is injective. Since this is a finite 
length  $Q(A)$-module, $b$ has to act invertibly. 
Hence 
\[ Q(A) \otimes_{A} H \cong Q(A) \otimes_{A} H \otimes_{B} Q(B) . 
\]
By symmetry we also get 
\[ H \otimes_{B} Q(B) \cong Q(A) \otimes_{A} H \otimes_{B} Q(B) . 
\]

By the Auslander property, for every $i>-d$ we have
\[ \opn{Cdim}_{R; A} \mrm{H}^{i} R 
=\opn{Cdim}_{R; A} \mrm{Ext}^{i}_A(A, R)
\leq -i < d= \opn{Cdim}_{R; A} A.\]
If $i<-d$, then the Auslander property shows that $\mrm{H}^{i} R=0$.
This implies that 
\[ Q(A)\otimes H^i R = H^i R \otimes Q(B) = 0 \]
for all $i \neq -d$.

\medskip \noindent 
(2) By the proof of part 1, $H = \mrm{H}^{-d} R$ 
is torsion-free on both sides and 
\[ Q(A)\otimes_A H \cong Q(A) \otimes_A H \otimes_B Q(B) \cong 
H\otimes_B Q(B) . \]
 
By \cite[2.15]{YZ1} the functor
$\opn{Ext}^{-d}_A(-, R)$ is a duality between the subquotient 
category $\msf{M}_{d, \mrm{f}} / \msf{M}_{d - 1, \mrm{f}}$ 
of $\cat{Mod} A$ and the subquotient 
category $\msf{M}^{\mrm{op}}_{d, \mrm{f}} / 
\msf{M}^{\mrm{op}}_{d - 1, \mrm{f}}$ 
of $\cat{Mod} B^{\mrm{op}}$. (For the definition of the 
category $\msf{M}_{d, \mrm{f}}$ see \cite[Section 2]{YZ1}.)
Now 
$\msf{M}_{d, \mrm{f}} / \msf{M}_{d - 1, \mrm{f}}
\approx \cat{Mod}_{\mrm{f}} Q(A)$, 
and similarly for $B^{\mrm{op}}$. We conclude that 
\[ \opn{Hom}_{Q(A)} 
\bigl(-, Q(A) \otimes_{A} H \otimes_B{Q(B)} \bigr) :
\cat{Mod}_{\mrm{f}} Q(A) \to
\cat{Mod}_{\mrm{f}} Q(B)^{\mrm{op}} \]
is a duality. Since for any torsion-free
$A$-module $M$ the Goldie rank $\opn{Grank}_A M$ coincides with 
the length of $Q(A) \otimes_A M$ in $\cat{Mod} Q(A)$, and likewise 
for $B^{\mrm{op}}$, we are finished.

\medskip \noindent
(3) The assertion follows from Lemma \ref{xx2.5}(2) and part (2).
\end{proof}


\section{The Trace Property}

To prove the results stated in the introduction we need to extend 
the formula in Proposition \ref{xx3.3}(3) to all prime ideals. 
First we need to construct a correspondence between the set of 
the prime ideals of $A$ and that of $B$. With this correspondence 
we show that the multiplicity $\mu(\fp)$ for $\fp \subset A$ can 
be computed by the Goldie rank of $B/\fq$ for the corresponding
prime $\fq \subset B$.

\begin{dfn} \label{xx4.1}
Suppose that $A \to \bar{A}$ and $B\to \bar{B}$ are algebra 
homomorphisms such that $\bar{A}$ is a finite left $A$-module and 
$\bar{B}$ is a finite right $B$-module. Also suppose $R$ 
and $\bar{R}$ are dualizing complexes over $(A, B)$ and 
$(\bar{A}, \bar{B})$ respectively. We say $R$ {\it induces} $\bar{R}$ 
if there exists a pair of morphisms 
\[ \opn{Tr}, \opn{Tr}^{\mrm{op}} : \bar{R} \to R \]
in $\msf{D}(\cat{Mod} A \otimes B^{\mrm{op}})$,
such that the induced morphisms
$\bar{R} \to \opn{RHom}_{A}(\bar{A}, R)$
in 
$\cat{D}(\cat{Mod} \bar{A} \otimes B^{\mrm{op}})$
and 
$\bar{R} \to \opn{RHom}_{B^{\mrm{op}}}(\bar{B}, R)$
in 
$\cat{D}(\cat{Mod} A \otimes \bar{B}^{\mrm{op}})$
are both isomorphisms. 
The morphism $\opn{Tr}$ is called a {\em left trace}, and 
$\opn{Tr}^{\mrm{op}}$ is called a {\em right trace}.
\end{dfn}

It might happen that 
$\opn{Tr} = \opn{Tr}^{\mrm{op}}$ 
(cf.\ the rigid trace in  \cite{YZ1}), but we do not require this.

The following two lemmas are proved in \cite[3.9]{YZ1}. 
(The rigidity mentioned in \cite[3.9]{YZ1} is not used in 
the proof.)

\begin{lem} \label{xx4.2} 
Suppose $R$ and $\bar{R}$ are dualizing complexes over $(A,B)$ 
and $(\bar{A}, \bar{B})$ respectively, and $\bar{R}$ is 
induced by $R$. Then for every $\bar{A}$-module $M$ one has
\[ \opn{RHom}_A(M, R) \cong \opn{RHom}_{\bar{A}}(M, \bar{R}) \]
in $\cat{D}(\cat{Mod} B^{\mrm{op}})$.
The same holds after exchanging $A$ and $B^{\mrm{op}}$.
\end{lem}

\begin{lem} \label{xx4.3}
Let $R$ be an Auslander dualizing complex over $(A, B)$. Suppose 
$\bar{R}$ is a dualizing complex over $(\bar{A}, \bar{B})$ induced 
by $R$. Then:
\begin{enumerate}
\item{} $\bar{R}$ is Auslander.
\item{} For every  $\bar{A}$-module $M$ one has
$\opn{Cdim}_{R; A} M = \opn{Cdim}_{\bar{R}; \bar{A}} M$; and
likewise after replacing $A$ with $B^{\mrm{op}}$.
\end{enumerate}
\end{lem}

\begin{dfn} \label{xx4.4}
A dualizing complex $R$ over $(A, B)$ is said to have the {\it 
trace property for ideals} if the following conditions hold:
\begin{enumerate}
\rmitem{i}
There is an inclusion-preserving bijection 
$\phi$ between the lattice of two-sided 
ideals of $A$ and the lattice of two-sided ideals of $B$.
\rmitem{ii} Given an ideal $\mfrak{a} \subset A$ let 
$\mfrak{b} := \phi(\mfrak{a})$, 
$\bar{A} := A / \mfrak{a}$ and $\bar{B} := B / \mfrak{b}$. 
Then there is a dualizing complex
$\bar{R}$ over $(\bar{A}, \bar{B})$ that is induced by $R$.
\end{enumerate}
If these conditions hold for the set of prime ideals (instead of
the set of all ideals) then $R$ is said to have the
{\it trace property for prime ideals}.
\end{dfn}

The trace property for ideals is not automatic (see Examples 
\ref{xx6.1} and \ref{xx6.2}). But this property can be proved for 
several classes of dualizing complexes. 

\begin{prop} \label{xx4.5} 
The following dualizing complexes have the trace 
properties for ideals.
\begin{enumerate}
\item{}
Suppose $A$ is a filtered algebra and $\opn{gr} A$ is noetherian 
connected with balanced dualizing complex. Then the rigid 
dualizing complex over $A$ has the trace property for ideals. The 
lattice isomorphism $\phi$ is the identity on the set of ideals 
of $A$.
\item{}
Let $A$ and $B$ be algebras satisfying the hypotheses in 
Lemma \ref{xx2.1}(3). Then every pre-balanced dualizing complex over 
$(A,B)$ has the trace property for ideals, and the lattice 
isomorphism $\phi$ is the one given in 
\cite[24.6]{AF} \tup{(}and in \cite[5.7]{CWZ}\tup{)}.
\item{}
Let $A$ be a complete semilocal algebra with $A / \fm$ finite 
over $\k$ for every maximal ideal $\fm$. Let $R$ be the balanced 
dualizing complex over $A$ \tup{(}see \cite[Definition 
3.7]{CWZ}\tup{)}. Then $R$ has the trace property for ideals. The 
lattice isomorphism $\phi$ is the identity on the set of ideals 
of $A$. 
\end{enumerate}
\end{prop}

\begin{proof} (1) This follows from \cite[6.17]{YZ1}.

\medskip \noindent
(2), (3) These assertions follow from \cite[5.6, 5.7 and 5.9]{CWZ}.
\end{proof}

To compute multiplicities of injectives
we only need the trace property for prime ideals.

\begin{lem} \label{xx4.6}
Let $R$ be an Auslander dualizing complex over $(A,B)$. If $R$ 
has the trace property for ideals then it has the trace property 
for prime ideals.
\end{lem}

\begin{proof} 
Suppose $\phi$ is a lattice isomorphism between
the set of all ideals of $A$ and the set of all ideals of $B$
as in Definition \ref{xx4.4}(1). We only need to show that the 
restriction of $\phi$ is a bijection between
the set of prime ideals of $A$ and the set of prime ideals of
$B$. By symmetry, it suffices to show that if $\mfrak{a} \subset A$
is prime then $\mfrak{b} := \phi(\mfrak{a}) \subset B$ is prime. 

Suppose $\mfrak{b}$ is not prime. Let $\bar{A} := A / \mfrak{a}$ and
$\bar{B} := B / \mfrak{b}$ be the quotient rings, and let
$\bar{R}$ be the induced dualizing complex. By Lemma \ref{xx4.3}
$\bar{R}$ is also Auslander. Suppose 
$\opn{Cdim}_{\bar{R}; \bar{A}} \bar{A} = q$. 
Then of course 
$\opn{Cdim}_{\bar{R}; \bar{B}^{\mrm{op}}} \bar{B} = q$. 
Since $\mfrak{b}$ is not prime, there is a 
prime ideal $\fq \supsetneqq \mfrak{b}$ such that 
$\opn{Cdim}_{\bar{R}; \bar{B}^{\mrm{op}}} \bar{B} / \fq = q$. 
Suppose $\fp \subset A$ is the ideal corresponding to $\fq$. 
Then $\fp \supsetneqq \mfrak{a}$ and 
$\opn{Cdim}_{\bar{R}; \bar{A}} \bar{A} / \fp = q$. 
This is impossible because $A / \mfrak{a}$ is prime of  
$\opn{Cdim}  = q$.
\end{proof}


\section{Multiplicity, Part II}

In this section we prove Theorems \ref{xx0.1} and \ref{xx0.2} and 
Corollaries \ref{xx0.3} and \ref{xx0.5}. Recall that 
$\opn{Grank} M$ denotes the Goldie rank of a module $M$.

\begin{thm} \label{xx5.1}
Let $R$ be an Auslander dualizing complex over a pair of rings
$(A, B)$. Assume $R$ has the trace property for prime ideals. 
Let $R \to I$ be a minimal injective resolution of $R$ 
as a complex of $A$-modules.
Let $\fp$ be a prime ideal of $A$, and 
let $\fq$ be the corresponding prime ideal of $B$. 
Denote by $\mu_i(\fp)$ be the multiplicity of the 
indecomposable injective module $J(\fp)$ in $I^i$.
Then 
\[ \mu_i(\fp) = \begin{cases} 
\opn{Grank} B / \fq & \text{\em if }  
i = -\opn{Cdim}_R A / \fp \\
0 & \text{\em otherwise} .
\end{cases} \]
\end{thm}

\begin{proof} 
By the trace property for prime ideals and Lemma \ref{xx4.3}
we reduce to the case when $\fp = 0$. Then the assertion follows 
from Proposition \ref{xx3.3}(3).
\end{proof}

Recall that the trace property holds for several
classes of algebras (see Proposition \ref{xx4.5}). 
Also note that the Auslander condition
is quite natural; cf.\ \cite{YZ1}.

\begin{proof}[Proof of Theorem \tup{\ref{xx0.1}}] 
This follows from Theorem \ref{xx5.1} and 
Proposition \ref{xx4.5}(1).
\end{proof}

\begin{proof}[Proof of Corollary \tup{\ref{xx0.2}}] 
If $\opn{gr} A$ is noetherian, and if it is FBN (or PI) or has 
enough normal elements, then $A$ has a rigid Auslander 
$\opn{GKdim}$-Macaulay dualizing complex $R$
\cite[6.9]{YZ1}. By definition, $\opn{GKdim}$-Macaulay means that 
$\opn{Cdim}_R = \opn{GKdim}$. 
Thus the assertion follows from Theorem \ref{xx0.1}.
\end{proof}

\begin{cor} \label{xx5.2}
Let $A$ be a filtered ring such that the associated graded ring 
$\opn{gr} A$ is connected graded, noetherian and Auslander 
Gorenstein. Let $\nu_i(\fp)$ denote the multiplicity of $J(\fp)$ 
in the $i$th term of a minimal injective resolution of $A$
as left $A$-module. Then 
\[ \nu_i(\fp) = \begin{cases}
\opn{Grank} A / \fp & \text{\em if } 
i = \opn{Cdim} A - \opn{Cdim} A / \fp , \\
0 & \text{\em otherwise} . 
\end{cases} \]
Here $\opn{Cdim} := \opn{Cdim}_R$, where $R$ is the rigid 
dualizing complex over $A$, which is Auslander.
\end{cor}

\begin{proof} 
By \cite[6.18]{YZ1} we know that $R \cong A^\sigma[n]$, where 
$\sigma$ is a filtered algebra 
automorphism of $A$. Then $n = \opn{Cdim} A$, and
$R \cong A[n]$ in $\msf{D}(\cat{Mod} A)$. 
Thus the assertion follows from Theorem \ref{xx0.1}.
\end{proof}

\begin{proof}[Proof of Corollary \tup{\ref{xx0.3}}] 
The graded ring $\opn{gr} U(\mfrak{g})$ is commutative. 
By \cite[6.9]{YZ1} one has 
$\opn{Cdim}= \opn{GKdim}$. 
It is also clear that 
$\opn{GKdim} U(\mfrak{g})= \opn{dim}_{\k} \mfrak{g}$.
Therefore the assertion follows from Corollary \ref{xx5.2}.
\end{proof}

Applying Theorem \ref{xx5.1} to the FBN rings we obtain the 
following result. Note that there might not exist a dualizing 
complex over FBN rings $(A,B)$ (even $A=B$). If it exists, 
the dualizing complex may not have the Auslander property and 
the trace property [Example \ref{xx6.1}]. 

\begin{cor} \label{xx5.3}
Let $A$ and $B$ be FBN rings, and let $R$ be an Auslander 
dualizing complex over $(A,B)$ with the trace property for prime 
ideals. Let $R \to I$ be a minimal injective resolution of $R$ as 
complex of $A$-modules. Then for every integer $i$ there is a 
decomposition
\[ I^i \cong \bigoplus J(\fp)^{\mu_i(\fp)} , \]
where the sum ranges over all prime ideals $\fp \subset A$ with 
$\opn{Cdim}_R  A / \fp = -i$. Furthermore 
\[ \mu_i(\fp) = \opn{Grank} B / \phi(\fp) < \infty, \] 
where $\phi$ is the lattice isomorphism from Definition 
\tup{\ref{xx4.4}}. In particular, for any such $\fp$ the 
indecomposable injective $J(\fp)$ appears only in the $i$th 
term $I^i$, where $i := -\opn{Cdim}_R A / \fp$.
\end{cor}

\begin{proof} This follows from Theorem 
\ref{xx5.1} and the one-to-one correspondence between
primes ideals and indecomposable injectives \cite[8.14]{GW}.
\end{proof}

Observe that Corollary \ref{xx5.3} is similar to Theorem 
\ref{xx0.4} with a slightly different hypothesis.

\begin{rem} \label{xx5.4} 
In Theorem \ref{xx0.4} and Corollary \ref{xx5.3} the $i$th term 
in the minimal injective resolution of the Auslander dualizing 
complex $R$ (as complex of left modules, and, by symmetry,
also as complex of right modules) is pure of 
$\opn{Cdim}_R = -i$. 
When this happens we say that $R$ has pure minimal 
injective resolutions. 
One of the main results in \cite{YZ2} is that in this case the 
Cousin complex $\mrm{E} R$ is a residual complex, and there is an
isomorphism $\mrm{E} R \cong R$ 
in the derived category 
$\msf{D}(\cat{Mod} A \otimes B^{\mrm{op}})$. 
When considered as a complex of left (or right)
modules, $\mrm{E} R$ is a minimal injective complex. Thus there 
is a decomposition
$(\mrm{E} R)^i \cong \bigoplus J(\fp)^{\mu_i(\fp)}$
as left modules, and a decomposition
$(\mrm{E} R)^i \cong \bigoplus J(\fq)^{\mu_i(\fq)}$
of right modules, where the sums are over the prime ideals 
$\fp \subset A$ and $\fq \subset B$ such that
$\opn{Cdim}_{R; A} A / \fp = 
\opn{Cdim}_{R; B^{\mrm{op}}} B / \fq = -i$.
\end{rem}

\begin{proof}[Proof of Corollary \tup{\ref{xx0.5}}] 
Let $d$ be the Krull dimension of $A$. By \cite[3.10]{SZ} 
the ring $A$ is Auslander Gorenstein, 
$\opn{Kdim}$-Macaulay of injective 
dimension $d$. Since $A$ is PI, the dualizing complex
$R := A[d]$ is bifinite (see Lemma \ref{xx2.1}). 
Then Theorem \ref{xx0.4} gives the form of $i$th 
term of the minimal injective resolution of $R$
as complex of left modules. Note that 
$\opn{Cdim}_R M = \opn{Kdim} M$
for any finite $A$-module $M$. It remains 
to show that multiplicities are correct.

The complex $R$ is a pre-balanced dualizing complex over $A$. By 
\cite[3.9]{CWZ}, $A$ admits a balanced dualizing complex, say 
$R'$. By Proposition \ref{xx4.5}(3), $R'$ 
has the trace property for ideals; so Theorem 
\ref{xx5.1} holds for $R'$. Now by Proposition 
\ref{xx2.3} we have
$R' \cong L \otimes_A R$ in $\msf{D}(\cat{Mod} A^{\mrm{e}})$
for some invertible bimodule $L$. Since $A$ is local the 
bimodule $L$ is isomorphic to $A^{\sigma}$ for some automorphism 
$\sigma$. Therefore 
$R' \cong R = A[d]$ in $\msf{D}(\cat{Mod} A)$, and we deduce that 
the resolution of $A$ has the correct multiplicities.
\end{proof}

Finally we briefly mention that the multiplicity of an injective 
not of the form $J(\fp)$ is usually zero or infinite. In the next 
proposition we show 
this for injective hulls of simple modules. 

\begin{prop} \label{xx5.5}
Let $R$ be a dualizing complex over $(A,B)$. Suppose that for 
every ideal $\mfrak{b} \subsetneqq B$ the quotient ring
$B / \mfrak{b}$ is not right artinian. 
Let $S$ be a simple $A$-module. If the injective hull
$E(S)$ appears in the $i$th 
term of the minimal injective resolution of $R$
as complex of $A$-modules, then its multiplicity is infinite.
\end{prop}

\begin{proof} 
Let $R \to I$ be the minimal injective resolution of $R$
as complex of $A$-modules. Since $S$ 
is simple the complex $\opn{Hom}_A(S, I)$ has zero differentials, 
and hence 
$\opn{Ext}^i_A(S, R) \cong \opn{Hom}_A(S, I^i)$. 
Define $M := \opn{Ext}^i_A(S, R)$. 
If the multiplicity of $E(S)$ in $I^i$ is a positive integer, say 
$n$, then 
\[ \opn{Hom}_A(S, I^i) \cong 
\opn{Hom}_A(S, S^{n}) \cong D^{n} , \]
where $D := \opn{Hom}_A(S, S)$ is a division ring. Therefore  
$M$ is a $(D^{\mrm{op}} \otimes B^{\mrm{op}})$-module which is 
finite over $B^{\mrm{op}}$ and finite length over 
$D^{\mrm{op}}$. By Lenagan's lemma \cite[7.10]{GW}
$M$ is a finite length $B^{\mrm{op}}$-module. 
Let $\mfrak{b} := \opn{Ann}_{B^{\mrm{op}}}(M)$. 
Then $B / \mfrak{b}$ is right artinian. This contradicts our
hypothesis. Therefore the multiplicity of $E(S)$ 
in $I^i$ has to be either $0$ or infinite.
\end{proof}


\section{Examples}

In the first example we show that not every dualizing complex 
is weakly bifinite.

\begin{exa} \label{xx6.1} 
Let $C \subset D$ be division rings such that 
$\opn{dim}_{C^{\mrm{op}}} D$ is finite 
and $\alpha := \opn{dim}_{C} D$ is infinite \cite[5.6.1]{Co}. 
Consider the upper triangular matrix ring 
$A := \begin{bmatrix} D &D \\ 0 & C \end{bmatrix}$. 
By \cite[Exc.\ 9 on p.\ 286]{AF} the ring $A$ is two-sided 
artinian (hence FBN), and by \cite[7.5.1]{MR} the global 
dimension of $A$ is $1$. So $R := A$ is a dualizing complex
over $A$. 

We claim that the dualizing complex $R$ is not Auslander. 
Let's identify $C$ and $D$ with the two simple quotient rings of 
$A$. Let $E(D)$ be the injective hull of the simple left 
$A$-module $D$. By \cite[Exercises 8 and 9 on p.\ 286]{AF}
the $A$-module $E(D)$ is not finite;
indeed, $E(D) / D \cong D \cong C^{(\alpha)}$, 
a direct sum of $\alpha$ copies of $C$. Pick an $A$-submodule 
$M \subset E(D)$ of length $l > 3$ 
(the length of $A$ is $3$). Then any 
$\phi \in \mrm{Hom}_{A}(M, A)$ can't be an injection. It follows 
that $D \subset \opn{Ker}(\phi)$. But then the induced 
homomorphism
$\bar{\phi} : M / D \cong C^{l - 1} \to A$
must be zero. We conclude that 
$\mrm{Hom}_{A}(M, A) = 0$ and hence
$\opn{Cdim}_{R; A} M = -1$. 
On the other hand $\opn{Hom}_A(D, A) \neq 0$, implying that
$\opn{Cdim}_{R; A} D = 0$. Since $D \subset M$ this 
shows $R$ {\it is not Auslander.}

The minimal injective resolution of $R$ as 
left $A$-module is:
\[  0 \to A \to E(D)^2 \to C^{(\alpha)} \to 0 . \]
 From this we see that
$\mrm{Hom}_{A}(C, R) = 0$ and 
\[ \mrm{Ext}^1_{A}(C, R) \cong \mrm{Hom}_{A}(C, C^{(\alpha)})
\cong C^{(\alpha)} \]
as $A^{\mrm{op}}$-modules. 
Hence {\it $R$ is not weakly bifinite}. Furthermore  
{\it $R$ does not have trace property for \tup{(}prime\tup{)}
ideals}.  
\end{exa}

Now suppose $A$ is a finite $\k$-algebra. Then $A^* := 
\opn{Hom}_{\k}(A, \k)$ is the rigid dualizing complex over $A$. 
It is clear that $A^*$ is Auslander, bifinite, 
$\opn{Cdim}$-symmetric, and 
has the trace property for ideals. But usually $A$ has many 
non-isomorphic dualizing complexes. The next examples, taken from 
\cite[5.4]{ASZ}, show that there exist algebras $A$ and
dualizing complexes $R$ over $A$ such that:
\begin{enumerate}
\item $R$ is not Auslander.
\item $R$ is Auslander, but it does not have the
trace property for prime ideals, nor is it $\opn{Cdim}$-weakly 
symmetric.
\end{enumerate}

\begin{exa} \label{xx6.2}
Let $A := \begin{bmatrix} \k & \k \\ 0 & \k \end{bmatrix}$. Then 
$A$ is Auslander regular of global dimension 1, so $R := A$ is an 
Auslander dualizing complex over $A$. 
Clearly $R$ is bifinite (see Lemma \ref{xx2.1}(1)). 
Let $\k_1$ and $\k_2$ be the simple quotient rings of $A$,
corresponding to the prime ideals 
$\fp_1 := \begin{bmatrix} 0 & \k \\ 0 & \k \end{bmatrix}$ 
and 
$\fp_2 := \begin{bmatrix} \k & \k \\ 0 & 0 \end{bmatrix}$ 
respectively. An easy computation shows that 
$\opn{Ext}^0_A(\k_1 , R) \neq 0$ and 
$\opn{Ext}^0_{A^{\mrm{op}}}(\k_1, R) = 0$. 
Hence {\it $R$ is not $\opn{Cdim}$-weakly symmetric} and 
{\it not pre-balanced.} 

Note that $\fp_1$, considered as left module, 
is an injective hull of the simple 
$A$-module $\k_1$. The minimal injective resolution of 
$A$ as left $A$-module is
\[ 0 \to A \to \fp_1 \oplus \fp_1  \to \k_2 \to 0 . \]
Hence  
\[ \opn{RHom}_A(\k_1, R) \cong
\opn{RHom}_A(\k_1, \k_1)^2 \cong \k_1 \oplus \k_1 \]
in $\msf{D}(\cat{Mod} \k_1 \otimes A^{\mrm{op}})$. 
Since every dualizing complex over $\k$ is isomorphic to 
$\k[n]$ for some integer $n$, it follows that 
$\opn{RHom}_A(\k_1, R)$ can not be a dualizing complex over 
$(\k_1, \k_1)$ nor over $(\k_1, \k_2)$. Therefore by Lemma 
\ref{xx4.2} $R$ {\it does not have trace property for prime ideals}. 
\end{exa} 

\begin{exa} \label{xx6.3}
Let $A := 
\begin{bmatrix} \k & \k^{n}\\ 0& \k \end{bmatrix}$
for some $n \geq 2$. Then $A$ has global dimension 1, and 
$R := A$ is a dualizing complex over $A$; but $R$ is 
not Auslander. Like in the previous example, 
$R$ is bifinite, and does not have the trace property for 
prime ideals.
\end{exa}

In \cite[5.4]{ASZ} the following was proved. Consider the ring $A$ 
of Example \ref{xx6.3}. Let $E$ be the injective hull of $A$ as 
$A^{\mrm{op}}$-module. Then $E$ has a ``natural'' 
$A$-module structure, but $E$ is not an injective $A$-module. 
We now prove a general statement,
which applies to Example \ref{xx6.1} as well.

\begin{prop} \label{xx6.4} 
Let $A$ be a noetherian ring of injective dimension $1$, which 
is not Auslander. Let $R := A$, a dualizing complex over 
$A$. Assume there is an exact dimension function 
$\opn{dim}$ on $\cat{Mod} A$ and 
$\cat{Mod} A ^{\mrm{op}}$ such that $R$ has 
pure minimal injective resolutions on both sides relative to 
$\opn{dim}$. Then there does not exist an $A$-bimodule 
$E$ that contains a sub-bimodule isomorphic to $A$, such that 
$E$ is an injective hull of $A$ on both sides.
\end{prop}

\begin{proof}
Suppose on the contrary that such a bimodule $E$ 
exists. Then the bimodule complex 
\[ R' := \bigl( 0 \to E \to E / A \to 0 \bigr),  \]
with $E$ in degree $0$,
is a minimal injective resolution of $R$ on both sides. 
By assumption $E$ and $E / A$ are $\opn{dim}$-pure on both sides. 
So $R'$ is a residual complex over $A$. By \cite[2.6]{Ye2}
the complex $R$ is Auslander (and $\opn{Cdim}_R = \opn{dim}$).
Hence $A$ is an Auslander ring, a contradiction.
\end{proof}

In the next example we discuss some properties of dualizing complexes
over Weyl algebras.

\begin{exa} \label{xx6.5}
Let $A$ be the $n$th Weyl algebra over $\k$, for some $n \geq 1$. 
Then $A$ is Auslander regular (and Cohen-Macaulay). Hence the 
bimodule $A$ is an Auslander (Cohen-Macaulay) dualizing complex over 
$A$. The rigid dualizing complex over $A$ is 
$A[2n]$; this is proved in \cite[2.6]{Ye4} for 
$\opn{char} \k = 0$, and for positive characteristic it follows 
from the fact that $A$ is Azumaya over its center and using 
\cite[Theorem 6.2]{YZ2}. 
By \cite{YZ3}, if $R$ is a dualizing complex over $(B, A)$ for
some left noetherian ring $B$, then $B$ is Morita equivalent to $A$
and $R \cong P[m]$ for some $m$ and some invertible bimodule $P$. 
We see that every dualizing complex over $(B, A)$ is Auslander 
(and Cohen-Macaulay). 

Now take $n = 1$. When $\opn{char} \k = 0$
the dualizing complex $A[1]$ is pre-balanced, and 
when $\opn{char} \k > 0$ the complex $A[2]$ is pre-balanced.
Thus every dualizing complex over $(B, A)$ is pre-balanced after
a suitable shift.

Suppose now $\opn{char} \k = 0$ and $n \geq 2$.
By a result of Stafford \cite{St} there is a simple 
$A$-module of $\opn{GKdim} > n$ (these are called non-holonomic 
simples). Using this fact we can show that the dualizing complex
$A[n]$ is not pre-balanced. This implies that $P[n]$ is not 
pre-balanced for any invertible bimodule $P$. For $m \neq n$ 
the dualizing complex $P[m]$ fails the pre-balanced test
for holonomic simple modules. Therefore there is no pre-balanced 
dualizing complex over $(B, A)$ for any noetherian ring $B$.

Finally, assume $\opn{char} \k = 0$ and $n \geq 1$. Let $S$ be a 
simple $A$-module. Then by Proposition \ref{xx5.5} the 
multiplicity of $E(S)$ in the $i$th term of the minimal injective 
resolution of $A$ as left module is either 0 or infinite.
\end{exa}


\section{Miscellaneous}

The results in this section are not directly related to  
multiplicities of injectives, but they are related to the 
methods and ideas presented in the previous sections. 

Recall that when we say $R$ is a dualizing complex over 
$(A,B)$ we assume that $A$ is left noetherian and $B$ is 
right noetherian. The first proposition is a consequence 
of Lemma \ref{xx2.5}.

\begin{prop} \label{xx7.1}
Let $R$ be a dualizing complex over a pair of rings $(A,B)$. Then 
$A$ is PI if and only if $B$ is PI.
\end{prop}

\begin{proof} 
Assume $B$ is PI; we want to show $A$ is PI too. Since $A$ is 
left noetherian, it suffices to show that $A / \fp$ is PI for all 
primes $\fp \subset A$. So choose one such $\fp$.
By \cite[1.11]{YZ1} the $i$th term in the
minimal injective resolution of $R$ as complex of $A$-modules 
contains a nonzero left ideal of $A / \fp$. By 
Lemma \ref{xx2.5}(3) we know that
$L := \opn{Ext}^i_A(A / \fp, R)$ is not torsion as 
$(A / \fp)$-module. Let $N$ be the torsion submodule of $L$, 
so that $L / N$ is a nonzero $(A / \fp)$-$B$-bimodule which 
is torsion-free as $(A / \fp)$-module. 
Pick a nonzero sub-bimodule $M \subset L / N$ which 
is a faithful $(B / \fq)^{\mrm{op}}$-module for some prime ideal 
$\fq \subset B$. Let $Q$ be the fraction ring of $B / \fq$. 
Since $B$ is PI the simple artinian ring $Q$ is finite over its 
center $\mrm{Z}(Q)$. Now  
$M' := M \otimes_{B/\fq} Q$ is a finite $Q^{\mrm{op}}$-module, 
which is also a torsion-free $(A / \fp)$-module. We obtain a ring 
injection 
$A / \fp \inj \opn{End}_Q(M')$.
Since the latter is a finite $\mrm{Z}(Q)$-algebra 
it follows that $A / \fp$ is PI.
\end{proof}

The second proposition is an application of the trace property.

\begin{prop} \label{xx7.2} 
Let $R$ be an Auslander dualizing complex over a pair of rings
$(A,B)$. Assume $R$ has the trace property for prime ideals. 
Then $A$ is left FBN if and only if $B$ is right FBN.
\end{prop}

\begin{proof} 
Suppose on the contrary that $A$ is left FBN but $B$ is not
right FBN. Then there is a prime ideal $\fq \subset B$ and a 
faithful finite $(B / \fq)^{\mrm{op}}$-module $M$ 
which is torsion \cite[8.2]{GW}. By the trace property for prime 
ideals we can assume that $A$ and $B$ are prime rings and
$\fq = 0$. Since $M$ is a torsion $B^{\mrm{op}}$-module we have
\[ \opn{Cdim}_{R; B^{\mrm{op}}} M < \opn{Cdim}_{R; B^{\mrm{op}}} B . 
\]
By the Auslander property we also know that
\[ \opn{Cdim}_{R; A} \opn{Ext}^i_{B^{\mrm{op}}}(M, R) <
\opn{Cdim}_{R; B^{\mrm{op}}} B = \opn{Cdim}_{R; A} A \]
for any $i$. 
Since $A$ is left FBN, there is a nonzero two-sided ideal 
$\mfrak{a} \subset A$ such that 
$\mfrak{a} \cdot \Ext^i_{B^{\mrm{op}}}(M, R)=0$
for all $i$.

Now let $N$ be a finite $A$-module which is not faithful. We claim 
that $\Ext^i_A(N, R)$ is a non-faithful $B^{\mrm{op}}$-module 
for all $i$. By noetherian induction and exact sequences we can 
reduce to the case when $\fp N = 0$ for some nonzero prime ideal 
$\fp \subset A$. By the trace property there is a nonzero prime ideal 
$\fq \subset B$ and 
an induced dualizing complex $\bar{R}$ over $(A / \fp, B / \fq)$.
Then 
\[ \Ext^i_A(N, R) \cong \Ext^i_{A / \fp}(N, \bar{R}) \]
as $B^{\mrm{op}}$-modules. It follows that
$\Ext^i_A(N, R) \cdot \fq = 0$.

Applying the above claim to $N := \Ext^j_{B^{\mrm{op}}}(M, R)$ 
one sees that the $B^{\mrm{op}}$-module 
$\Ext^i_A(\Ext^j_{B^{\mrm{op}}}(M, R), R)$ is not faithful for 
all $i$ and $j$. Finally the double-Ext spectral sequence 
\cite[1.7]{YZ1} shows that $M$ is not faithful, a contradiction.
\end{proof}

If $A=B$ we have the following corollary.

\begin{cor} \label{xx7.3}
Let $A$ be a noetherian algebra. Suppose $A$ admits an Auslander 
dualizing complex with the trace property for prime ideals. Then 
$A$ is left FBN if and only if it is right FBN.
\end{cor}

It is unknown whether every noetherian left FBN ring is also right
FBN \cite[p.\ 132]{GW}.

Finally we show that the ``$\opn{Cdim}$-symmetric'' property is 
sometimes a consequence of other properties. The following 
proposition is an analogue of \cite[3.10]{SZ}. The proof of this 
is also very similar to the proof of \cite[3.10]{SZ}. For the 
convenience of reader we include it here.

\begin{prop} \label{xx7.4}
Let $A$ and $B$ be FBN rings and let $R$ be a pre-balanced,
weakly bifinite dualizing complex over $(A, B)$. Then:
\begin{enumerate}
\item $R$ is Auslander.
\item $\opn{Cdim}_{R; A} M = \opn{Kdim} M$ for all finite
$A$-modules $M$, and likewise for all finite 
$B^{\mrm{op}}$-modules. As a consequence, $R$ is
$\opn{Cdim}$-symmetric. 
\end{enumerate}
\end{prop}

\begin{proof} 
We will prove the following four statements by induction on the
Krull dimension of the finite $A$-module $M$.
\begin{enumerate}
\rmitem{a} $\Ext^i_A(M, R) = 0$ if $i > 0$ or $i < -\opn{Kdim} M$.
\rmitem{b} $\opn{Cdim}_{R; A} M = \opn{Kdim} M$.
\rmitem{c} $\opn{Kdim} \Ext^j_A(M, R) = \opn{Kdim} M$ where 
$j := j_R(M)$.
\rmitem{d} For every $i \geq j_R(M)$  one has 
$\opn{Kdim} \Ext^i_A(M, R) \leq -i$.
\end{enumerate}
By symmetry these statements will also hold after replacing $A$ 
with $B^{\mrm{op}}$. The Auslander condition will then follow.

If $\opn{Kdim} M=0$ then $M$ has finite length, and the problem
can be reduced to the case when $M$ is simple, in which 
case the assertions follow from the pre-balanced hypothesis.

Now suppose $\opn{Kdim} M = d > 0$, and assume the assertions hold
for all modules of lower Krull dimension. Using noetherian 
induction and \cite[Lemma 2.1]{SZ}, it suffices to check the case
$M = A / \fp$ for some prime ideal $\fp$. Let 
$\bar{A}$ be the quotient ring $A / \fp$. For every
regular element $u \in \bar{A}$ right multiplication gives a short
exact sequence
\[ 0 \to \bar{A} \to \bar{A} \to \bar{A} / \bar{A} u \to 0 , \]
which induces a long exact sequence
\begin{equation} \label{eqnx}
\Ext^i_A(\bar{A} / \bar{A} u, R) \to 
\Ext^i_A(\bar{A}, R) \to \Ext^i_A(\bar{A}, R) \to 
\Ext^{i+1}_A(\bar{A} / \bar{A} u , R) . 
\end{equation}
Observe that 
$\opn{Kdim} \bar{A} / \bar{A} u < d$. 
Write $E = \Ext^i_A(\bar{A}, R)$ and let 
$Q$ be the fraction ring of $\bar{A}$.

Suppose $i > 0$ or $i < -d$. In this case the two end terms 
in (\ref{eqnx}) are zero by the induction hypothesis, so
left multiplication by $u$ on $E$ is an isomorphism. 
Since $Q$ is obtained by inverting all such
regular elements $u$ it follows that $E$ is a $Q$-module. 
By the ``weakly bifinite'' hypothesis $E$ is a finite 
$\bar{A}$-module. If $E \neq 0$ then
$Q$ is a submodule of a finite direct sum of copies of $E$. This 
implies that $Q$ is a finite $\bar{A}$-module, and hence
$Q = \bar{A} v^{-1}$ for some regular element $v \in \bar{A}$. 
Then $Q = Q v = \bar{A} v^{-1} v = \bar{A}$, 
so $\fp$ is a maximal ideal $\opn{Kdim} \bar{A} = 0$, 
which is a contradiction. Therefore $E = 0$; this is statement
(a). 

Now suppose $-d \leq i \leq 0$. By the inductive hypothesis 
\[ \opn{Kdim} \Ext^{i+1}_A (\bar{A} / \bar{A} u, R)
\leq -i - 1 . \]
The long exact sequence \ref{eqnx} shows that 
$E / u E$ is a submodule of 
$\Ext^{i+1}_A(\bar{A} / \bar{A} u, R)$, so 
\[ \opn{Cdim}_{R; B^{\mrm{op}}} E / u E \leq -i - 1 . \]
From the induction hypothesis (this time $B^{\mrm{op}}$-modules) 
we see that
$\opn{Kdim} E / u E \leq -i - 1$. 
Since this holds for all regular elements $u$, the argument in 
\cite[p.\ 1003]{SZ} shows that 
$\opn{Kdim} E \leq -i$. 
This is statement (d). 

Finally we look at the case $i = -d$. The last paragraph shows that 
$\opn{Kdim} E \leq d$. We want to show there is equality. If not 
then $\opn{Kdim} E < d$, and by the previous two paragraphs
$\opn{Kdim} \Ext^{q}_A(\bar{A},R) < d$ for all $q$. 
By the induction hypothesis for $B^{\mrm{op}}$-modules we see that 
$\Ext^p_{B^{\mrm{op}}}(\Ext^q_A(\bar{A}, R), R)$ 
has Krull dimension less than $d$ for all $p,q$. 
The double-Ext spectral sequence 
\cite[1.7]{YZ2} implies that $\opn{Kdim} \bar{A} < d$,
and this is a contradiction. Therefore 
$\opn{Kdim} \Ext^{-d}_A(\bar{A}, R) = d$,
and as a consequence we obtain statements (b) and (c). 
\end{proof}

%


\end{document}